\documentclass[journal]{IEEEtran}
\usepackage[dvipdfmx]{graphicx,xcolor}
\usepackage{framed}
\usepackage{bm}
\usepackage{comment}
\usepackage{cite}
\usepackage{algorithm}
\usepackage{algorithmic}
\usepackage{slashbox}

\newtheorem{remark}{Remark}

\usepackage{mathtools}
\mathtoolsset{showonlyrefs=true}

\usepackage{ascmac}

\ifCLASSINFOpdf
\else
\fi
\usepackage{amsmath}
\begin{document}
%
\title{Structure-preserving $H^2$ optimal model reduction\\ based on Riemannian trust-region method}
%
%
%

\author{Kazuhiro Sato and Hiroyuki Sato
\thanks{K. Sato is with the School of Regional Innovation and Social Design Engineering,
 Kitami Institute of Technology,
 Hokkaido 090-8507, Japan,
email: ksato@mail.kitami-it.ac.jp}
\thanks{H. Sato is with the Department of Information and Computer Technology, Tokyo University of Science, Tokyo, 125-8585 Japan,
email: hsato@rs.tus.ac.jp} 
}

%
%

\markboth{This paper was published in IEEE Transactions on Automatic Control (DOI: 10.1109/TAC.2017.2723259)}%
{Shell \MakeLowercase{\textit{et al.}}: Bare Demo of IEEEtran.cls for IEEE Journals}
%



\maketitle

\begin{abstract}
This paper studies stability and symmetry preserving $H^2$ optimal model reduction problems of linear systems which include  linear gradient systems as a special case.
The problem is formulated as a nonlinear optimization problem on the product manifold of the manifold of symmetric positive definite matrices and the Euclidean spaces.
To solve the problem by using the trust-region method,
the gradient and Hessian of the objective function are derived.
Furthermore, it is shown that if we restrict our systems to gradient systems, the gradient and Hessian can be obtained more efficiently.
More concretely, by symmetry, we can reduce linear matrix equations to be solved. In addition, by a simple example, we show that 
the solutions to our problem and a similar problem in some literatures are not unique and the solution sets of both problems do not contain each other in general.
Also, it is revealed that the attained optimal values do not coincide.
Numerical experiments show that the proposed method gives a reduced system with the same structure with the original system although the balanced truncation method does not.
\end{abstract}


%
\IEEEpeerreviewmaketitle

\section{Introduction}
%
%
%
%
\IEEEPARstart{M}{odel} reduction method reduces the dimension of the state of a given system to facilitate the controller design.
The most famous method is called balanced truncation method, which gives a stable reduced order model with guaranteed $H^{\infty}$ error bounds \cite{antoulas2005approximation, moore1981principal}.
Another famous method is moment matching method \cite{astolfi2010model, ionescu2014families}, which gives a reduced system matching some coefficients of the transfer function of a given linear system.
In particular, \cite{ionescu2013moment, scherpen2011balanced} have discussed structure preserving model reduction methods which give a reduced model with the same structure with the original system.
However, the previous methods  do not guarantee any optimality.

In \cite{sato2015riemannian, yan1999approximate}, $H^2$ optimal model reduction problems have been studied.
Reference \cite{yan1999approximate} has formulated the problem as an optimization problem for minimizing the $H^2$ norm performance index subject to orthogonality constraints.
The constrained optimization problem can be regarded as an unconstrained problem on the Stiefel manifold.
To solve the problem, an iterative gradient flow method has been proposed in \cite{yan1999approximate}.
Reference \cite{sato2015riemannian} has reformulated the optimization problem on the Stiefel manifold into that on the Grassmann manifold
because the objective function of the $H^2$ optimal model reduction problem is invariant under actions of the orthogonal group.
To solve the problem, the Riemannian trust-region methods have been proposed in \cite{sato2015riemannian}.
The proposed methods in \cite{sato2015riemannian, yan1999approximate} preserve stability and symmetric properties of the original system.
However, the methods may not give good reduced models as shown in Section \ref{sec4} in this paper.

To obtain better reduced models,
this paper further exploits the stability and symmetry preserving $H^2$ optimal model reduction problem of linear systems which include gradient systems as a special case.
In particular, the present paper reformulates the problems in \cite{sato2015riemannian, yan1999approximate}  as an optimization problem on the product manifold of the manifold of  symmetric positive definite matrices and two Euclidean spaces.
This novel approach is the first to the best of our knowledge.
A global optimal solution in this formulation gives a smaller value of the objective function than that in  \cite{sato2015riemannian, yan1999approximate}.
Furthermore, the $H^2$ optimal model reduction problem of gradient systems is formulated as another specific optimization problem.

The contributions of this paper are as follows.\\
1)  We derive the gradient and Hessian of the objective function of the new optimization problem.
By using them, we can apply the Riemannian trust-region method to solve the new problem.
Furthermore, it is shown that if we restrict our systems to the gradient systems, the gradient and Hessian can be obtained more efficiently.
More concretely, by symmetry, we can reduce linear matrix equations to be solved.
Some numerical experiments demonstrate that the proposed Riemannian trust-region method gives a reduced system which is sufficiently close to the original system,
even if the balanced truncation method and the method in \cite{sato2015riemannian} do not. \\
2) By a simple example, we show that
the solutions to our problem and the problem in \cite{sato2015riemannian, yan1999approximate} are not unique and the solution sets of both problems do not in general contain each other.
Also, it is revealed that the attained optimal values do not coincide in general.

This paper is organized as follows.
In Section \ref{sec2}, we formulate the structure preserving $H^2$ optimal model reduction problem on the manifold.
In Section \ref{sec3}, we first review the geometry of the manifold of the symmetric positive definite matrices.
Next, we derive the Euclidean gradient of the objective function and then we give the Riemannian gradient and Riemannian Hessian to develop the trust-region method.
In Section \ref{sec_Pro2}, the $H^2$ optimal model reduction problem of gradient systems are discussed.
In Section \ref{secnew_5}, we study the difference between our problem and the problem in \cite{sato2015riemannian, yan1999approximate}.
Section \ref{sec4} shows some numerical experiments to investigate the performance of the proposed method.
We demonstrate that the objective function in the case of the proposed method takes a smaller value than those of the balanced truncation method and the method in \cite{sato2015riemannian}.
Furthermore, the experiments indicate that although the balanced truncation method does not preserve the original structure,
the proposed method does.
The conclusion is presented in Section \ref{sec5}.

{\it Notation:} The sets of real and complex numbers are denoted by ${\bf R}$ and ${\bf C}$, respectively.
The identity matrix of size $n$ is denoted by $I_n$.
The symbols ${\rm Sym}(n)$ and ${\rm Skew}(n)$ denote the sets of symmetric and skew-symmetric matrices in ${\bf R}^{n\times n}$, respectively.
The set of  symmetric positive definite matrices in ${\bf R}^{n\times n}$ is denoted by ${\rm Sym}_+(n)$.
The symbols $GL(r)$ and $O(r)$ are the general linear group and the orthogonal group of degree $r$, respectively.
Given a matrix $A\in {\bf R}^{n\times n}$, ${\rm tr} (A)$ denotes the sum of the diagonal elements of $A$
and ${\rm sym}(A)$ denotes the symmetric part of $A$; i.e., ${\rm sym} (A) = \frac{A+A^T}{2}$.
Here, $A^T$ denotes the transposition of $A$.
The tangent space at $x$ on a manifold $\mathcal{M}$ is denoted by $T_x \mathcal{M}$.
Given a smooth function $f$ on a manifold $\mathcal{M} \subset {\bf R}^{n_1\times n_2}$, the symbol $\bar{f}$ is the extension of $f$ to the ambient Euclidean space ${\bf R}^{n_1\times n_2}$.
The symbols $\nabla$ and ${\rm grad}$ denote the Euclidean and Riemannian gradients, respectively; i.e., given a smooth function $f$ on a manifold $\mathcal{M} \subset {\bf R}^{n_1\times n_2}$,
$\nabla$ and  ${\rm grad}$ act on  ${\bar f}$ and $f$, respectively.
The symbol ${\rm Hess}$ denotes the Riemannian Hessian.
Given a transfer function $G$, $||G||_{H^2}$ denotes the $H^2$  norm of $G$.

\section{Problem setup} \label{sec2}

We consider the $H^2$ optimal model reduction problem of a linear time invariant system
\begin{align}
\begin{cases}
\dot{x} = -Ax +Bu, \\
y = Cx,
\end{cases} \label{1}
\end{align}
where $x\in {\bf R}^n$, $u\in {\bf R}^m$, and $y\in {\bf R}^p$ are the state, input, and output, respectively,
and where $A\in {\bf R}^{n\times n}$, $B\in {\bf R}^{n\times m}$, and $C\in {\bf R}^{p\times n}$ are constant matrices.
Throughout this paper, we assume $A\in {\rm Sym}_{+}(n)$, and thus, all the eigenvalues of $-A$ are negative.
Thus, the original system has a stable and symmetric state transition matrix. 
Note that if $m=p$ (i.e., the number of the output variables is the same with that of the input variables) and $B^T=C$, then the system \eqref{1} is a linear gradient system \cite{ionescu2013moment, scherpen2011balanced}.
Note also that the following discussion fully exploits the symmetry of $A$ and does not apply to systems with a non-symmetric matrix $A$.

The structure preserving $H^2$ optimal model reduction problem in this paper is to find $A_r\in {\rm Sym}_+(r)$, $B_r\in {\bf R}^{r\times m}$, and $C_r\in {\bf R}^{p\times r}$ for
a fixed integer $r$ $(<n)$ such that the associated reduced system
\begin{align}
\begin{cases}
\dot{x}_r = -A_rx_r +B_ru, \\
y_r = C_rx_r
\end{cases} \label{2}
\end{align}
best approximates the original system \eqref{1} in the sense that the $H^2$ norm of the transfer function of the error system between the original system \eqref{1} and the reduced system \eqref{2} is minimized.
That is, the stability and symmetry of the state transition matrix are preserved because
the reduced matrices $A_r$, $B_r$, and $C_r$ have the same structures with the  original matrices $A$, $B$, and $C$, respectively.
Note that the symmetry preservation is significant because the symmetry implies that any oscillations never occur when $u=0$.
This is because all the eigenvalues of any symmetric matrices are real numbers. 
If the state transition matrix of the reduced system \eqref{2} is not symmetric, some oscillations may be observed under $u=0$ in contrast to the case of the original system \eqref{1}.

The optimization problem to be solved is stated as follows.
\begin{framed}
Problem 1:
\begin{align*}
&{\rm minimize} \quad J(A_r,B_r,C_r), \\
&{\rm subject\, to} \quad (A_r,B_r,C_r)\in M.
\end{align*}
\end{framed}

\noindent
Here, 
\begin{align}
\label{eq_J}
J(A_r,B_r,C_r):=||G-G_r||_{H^2}^2,
\end{align}
where
\begin{align*}
G(s) := C(sI_n +A)^{-1}B, \quad s\in {\bf C}
\end{align*}
is the transfer function of the original system \eqref{1} and
 $G_r$ is the transfer function of the reduced system \eqref{2}, and
\begin{align*}
M:= {\rm Sym}_+(r) \times {\bf R}^{r\times m} \times {\bf R}^{p\times r}.
\end{align*}

Since all the eigenvalues of $-A$ and $-A_r$ are negative,
the objective function $J(A_r,B_r,C_r)$ can be expressed
 as
\begin{align*}
J(A_r,B_r,C_r) &= {\rm tr}( C\Sigma_c C^T +C_rPC_r^T-2C_rX^TC^T) \\
&= {\rm tr} ( B^T\Sigma_o B +B_r^TQB_r + 2B^TYB_r),
\end{align*}
where the matrices $\Sigma_c$, $\Sigma_o$, $P$, $Q$, $X$, and $Y$ are the solutions to
\begin{align}
A\Sigma_c +\Sigma_cA -BB^T  &= 0,  \nonumber\\
A\Sigma_o + \Sigma_o A- C^TC &= 0, \nonumber\\
A_rP+PA_r-B_rB_r^T &=0, \label{3} \\
A_rQ+QA_r-C_r^TC_r &=0, \label{4} \\
AX+XA_r-BB_r^T &=0, \label{5} \\
AY+YA_r+C^TC_r &=0, \label{6}
\end{align} 
respectively.  A similar discussion can be found in \cite{sato2015riemannian}, which contains a more detailed explanation of the calculation.

As mentioned earlier, if $m=p$ and $B^T=C$, then the system \eqref{1} is a stable gradient system \cite{ionescu2013moment, scherpen2011balanced}.
If this is the case, Problem 1 can be replaced with the following problem.
\begin{framed}
Problem 2:
\begin{align*}
&{\rm minimize} \quad \tilde{J}(A_r,B_r), \\
&{\rm subject\, to} \quad (A_r,B_r)\in \tilde{M}.
\end{align*}
\end{framed}

\noindent
Here, 
\begin{align*}
\tilde{J}(A_r,B_r) := || B^T(sI_n+A)^{-1}B - B_r^T(sI_r+A_r)B_r||^2_{H^2}
\end{align*}
and 
\begin{align*}
\tilde{M}:= {\rm Sym}_+(r)\times {\bf R}^{r\times m}.
\end{align*}

We develop optimization algorithms for solving Problems 1 and 2 in Sections \ref{sec3} and \ref{sec_Pro2}, respectively.

\begin{remark}
We can also consider the reduced system expressed by
\begin{align*}
\begin{cases}
\dot{x}_r = -U^TAUx_r +U^TB u, \\
y_r = CU x
\end{cases}
\end{align*}
for $U$ belonging to the Stiefel manifold ${\rm St}(r,n):= \left\{ U\in {\bf R}^{n\times r}\, |\, U^TU =I_r\right\}$.
Then, Problem 1 is replaced with the following optimization problem on the Stiefel manifold.
\begin{framed}
Problem 3:
\begin{align*}
&{\rm minimize} \quad J(U^TAU,U^TB,CU), \\
&{\rm subject\, to} \quad  U \in {\rm St}(r,n).
\end{align*}
\end{framed}

\noindent
Reference \cite{sato2015riemannian} has proposed the trust-region method for solving Problem 3.
\end{remark}

\begin{remark}
It is beneficial to consider Problem 1 instead of Problem 3.
In fact,
if $U_*$ is a global optimal solution to Problem 3, then $(U^T_*AU_*,U_*^TB,CU_*)\in M$ is a feasible solution to Problem 1; i.e.,
the minimum value of Problem 3 is not smaller than that of Problem 1.
In Section \ref{sec4}, we verify this fact numerically.
Furthermore, we give an example in Section \ref{secnew_5} which shows that the critical points of Problems 1 and 3 do not necessarily coincide with each other.
\end{remark}

\begin{remark}
A possible drawback of Problem 1 is that the problem may not have a solution because the Riemannian manifold $M$ is not compact.
However, we always obtained solutions by using the trust-region method in this paper.
We leave a general mathematical analysis of the existence of the solution to Problem 1 to future work. 
\end{remark}


\section{Optimization algorithm for Problem 1} \label{sec3}

\subsection{General Riemannian trust-region method}
We first review Riemannian optimization methods including the Riemannian trust-region method following \cite{absil2009optimization} for readability of the subsequent subsections.
We also refer to \cite{absil2009optimization} for schematic figures of Riemannian optimization.
In this subsection, we consider a general Riemannian optimization problem to minimize an objective function $h$ defined on a Riemannian manifold $\mathcal{M}$.

In optimization on the Euclidean space $\mathcal{E}$, we can compute a point $x_+ \in \mathcal{E}$ from the current point $x \in \mathcal{E}$ and the search direction $d \in \mathcal{E}$ as $x_+=x+d$.
However, this update formula cannot be used on $\mathcal{M}$ since $\mathcal{M}$ is not generally  a Euclidean space.
For $x \in \mathcal{M}$ and $\xi \in T_x \mathcal{M}$, $x+\xi$ is not defined in general.
Even if $\mathcal{M}$ is a submanifold of the Euclidean space $\mathcal{E}$ and $x+\xi$ is defined as a point in $\mathcal{E}$, it is not generally on $\mathcal{M}$. 
Therefore, we seek for a next point $x_+$ on a curve on $\mathcal{M}$ emanating from $x$ in the direction of $\xi$.
Such a curve is defined by using a map called an exponential mapping ${\rm Exp}$, which is defined by a curve called geodesic.
More concretely, for any $x, y \in \mathcal{M}$ on a geodesic which are sufficiently close to each other, the piece of geodesic between $x$ and $y$ is the shortest among all curves connecting the two points.
For any $\xi\in T_x\mathcal{M}$, there exists an interval $I \subset {\bf R}$ around $0$ and a unique geodesic $\Gamma_{(x,\xi)}: I\rightarrow \mathcal{M}$ such that $\Gamma_{(x,\xi)}(0)=x$ and $\dot{\Gamma}_{(x,\xi)}=\xi$.
The exponential mapping ${\rm Exp}$ at $x\in \mathcal{M}$ is then defined through this curve as
\begin{align}
{\rm Exp}_x(\xi):=\Gamma_{(x,\xi)}(1). \label{exp_def}
\end{align}
This definition is well-defined because the geodesic $\Gamma_{(x,\xi)}$ has the homogeneity property $\Gamma_{(x,a\xi)}(t)=\Gamma_{(x,\xi)}(at)$ for any $a\in {\bf R}$ satisfying $at \in I$.
We can thus compute a point $x_+$ as
\begin{equation}
x_+= {\rm Exp}_x(\xi). \label{next_step}
\end{equation}

In the trust-region method, at the current point $x$, we compute a second-order approximation of the objective function $h$ based on the Taylor expansion.
We minimize the second-order approximation in a ball of a radius called trust-region radius.
See Section \ref{Sec:3.E} for detail discussion on the trust-region method for our problem.

We thus need the first and second-order derivatives of $h$, which are characterized by the Riemannian gradient and Hessian of $h$.
Since $\mathcal{M}$ is a Riemannian manifold, $\mathcal{M}$ has a Riemannian metric  $\langle \cdot, \cdot\rangle$,
which endows the tangent space $T_x \mathcal{M}$ at each point $x \in \mathcal{M}$ with an inner product $\langle \cdot, \cdot \rangle_x$.
The gradient ${\rm grad}\, h(x)$ of $h$ at $x \in \mathcal{M}$ is defined as a tangent vector at $x$ which satisfies
\begin{equation}
{\rm D}h(x)[\xi] = \langle {\rm grad}\, h(x), \xi\rangle_x \label{grad_def}
\end{equation}
for any $\xi \in T_x \mathcal{M}$.
Here, the left-hand side of \eqref{grad_def} denotes the directional derivative of $h$ at $x$ in the direction $\xi$.
The Hessian ${\rm Hess}\, h(x)$ of $h$ at $x$ is defined via the covariant derivative of the gradient ${\rm grad}\, h(x)$. If $\mathcal{M}$ is a Riemannian submanifold of Euclidean space, we can compute the Hessian ${\rm Hess}\, h(x)$ by using the gradient ${\rm grad}\, h(x)$ and the orthogonal projection onto the tangent space $T_x \mathcal{M}$.
Based on this fact, we derive the gradient and Hessian of our objective function $J$ in Section \ref{Sec:3.D}.

\subsection{Difficulties when we apply the general Riemannian trust-region method}

This subsection points out difficulties when we apply the above general Riemannian trust-region method.

The first difficulty is to obtain the geodesic $\Gamma_{(x,\xi)}$.
In fact, to get $\Gamma_{(x,\xi)}$, we may need to solve a nonlinear differential equation in a local coordinate system around $x\in \mathcal{M}$.
The  equation may only be approximately solved by a numerical integration scheme.
The numerical integration consumes a large amount of time in many cases.
As a result, it is difficult to obtain the exponential map defined by \eqref{exp_def} in general.

The second difficulty is how to choose 
a Riemmanian metric $\langle \cdot, \cdot \rangle$.
Since the gradient ${\rm grad}\,h(x)$ defined by \eqref{grad_def} varies by the Riemannian metric,
we should adopt a metric in such a manner that we can obtain the gradient in a short time. 
However, the adoption may imply that the manifold $\mathcal{M}$ is not geodesically complete.
Here, a Riemannian manifold is called geodesically complete if the exponential mapping is defined for every tangent vector at any point.
If $\mathcal{M}$ is not geodesically complete, we have to carefully choose $\xi$ in \eqref{next_step} in such a manner that $x_+$ is contained in $\mathcal{M}$.
This leads to computational inefficiency.

For example, consider the manifold ${\rm Sym}_{+}(r)$.
Since ${\rm Sym}_{+}(r)$ is a submanifold of the vector space ${\rm Sym}(r)$,
we can consider the induced metric $\langle \cdot,\cdot \rangle$ from the natural inner product in the ambient space ${\rm Sym}(r)$ as
\begin{align}
\langle \xi_1, \xi_2 \rangle_{S}:= {\rm tr}(\xi_1 \xi_2) \label{normal_metric}
\end{align}
for $\xi_1, \xi_2\in T_S {\rm Sym}_+(r)$. Here, $T_S {\rm Sym}_+(r)\cong {\rm Sym}(r)$ as explained in \cite{lang1999fundamentals}.
Then, the exponential map is simply given by
${\rm Exp}_S(\xi) = S+\xi.$
However, $S+\xi \not \in {\rm Sym}_+(r)$ for some $\xi\in T_S {\rm Sym}_+(r)$ because of lack of positive definiteness.
This means that ${\rm Sym}_+(r)$ is not geodesically complete.
As a result, we have to carefully choose $\xi \in T_S {\rm Sym}_+(r)$.

Our following discussion overcomes these difficulties.

\subsection{Geometry of the manifold ${\rm Sym}_{+}(r)$} \label{Sec:3.B}

This subsection introduces another Riemannian metric on the manifold ${\rm Sym}_{+}(r)$ \cite{lang1999fundamentals, Gallier2016, helgason1979differential, helmke1996optimization,  pennec2006riemannian}.
This is useful to develop an optimization algorithm for solving Problem 1 for the following reasons:\\
1) The geodesic is given by a closed-form expression. That is, we do not have to integrate a nonlinear differential equation.\\
2) The manifold ${\rm Sym}_+(r)$ is then geodesically complete in contrast to the case of  the Riemannian metric \eqref{normal_metric}. 
That is, ${\rm Exp}_S(\xi)\in {\rm Sym}_+(r)$ is always defined for any $\xi\in T_S{\rm Sym}_+(r)$.

For $\xi_1$, $\xi_2\in T_S {\rm Sym}_{+}(r)$, we define the Riemannian metric as
\begin{align}
\langle \xi_1, \xi_2\rangle_S := {\rm tr} ( S^{-1} \xi_1 S^{-1} \xi_2  ), \label{metric}
\end{align}
which is invariant under the group action $\phi_g: S \to gSg^T$ for $g \in GL(r)$; i.e.,
$\langle {\rm D} \phi_g(\xi_1), {\rm D} \phi_g(\xi_2) \rangle_{\phi_g(S)} = \langle \xi_1,\xi_2 \rangle_S$,
where the map ${\rm D}\phi_g: T_S {\rm Sym}_+(r) \rightarrow T_S {\rm Sym}_+(r)$ is a derivative map given by ${\rm D}\phi_g(\xi) =g\xi g^T$.
The proof that \eqref{metric} is a Riemannian metric can be found in Chapter XII in \cite{lang1999fundamentals}.

Let $f: {\rm Sym}_+(r) \rightarrow {\bf R}$ be a smooth function and $\bar{f}$ the extension of $f$ to the Euclidean space ${\bf R}^{r\times r}$.
The relation of the Euclidean gradient $\nabla \bar{f}(S)$ and the directional derivative ${\rm D}\bar{f}(S)[\xi]$ of $\bar{f}$ at $S$ in the direction $\xi$ is given by
\begin{align}
{\rm tr}\, (\xi^T\nabla \bar{f}(S)) ={\rm D} \bar{f}(S)[\xi].
\end{align}

\noindent
The Riemannian gradient ${\rm grad}\, f(S)$ is given by
\begin{align}
\langle {\rm grad} f (S), \xi \rangle_S &= {\rm D}f(S)[\xi] \\
&= {\rm D}\bar{f}(S)[\xi] \\
&= {\rm tr}\, (\xi^T {\rm sym}(\nabla \bar{f}(S)) ). \label{relation}
\end{align}
Here, we have used $\xi=\xi^T$.
From \eqref{metric} and \eqref{relation}, we obtain
\begin{align}
{\rm grad}\, f(S) = S {\rm sym}(\nabla \bar{f}(S)) S. \label{gradient}
\end{align}

According to Section 4.1.4 in \cite{jeuris2012survey},
the Riemannian Hessian ${\rm Hess}\,f(S): T_S {\rm Sym}_+(r) \rightarrow T_S {\rm Sym}_+(r)$ of the function $f$ at $S\in {\rm Sym}_+(r)$ is given by
\begin{align}
{\rm Hess}\,f(S)[\xi] = {\rm D}{\rm grad}\, f(S)[\xi] - {\rm sym} ({\rm grad}\, f(S) S^{-1} \xi ). \label{Hess1}
\end{align}
Hence, \eqref{gradient} and \eqref{Hess1} yield
\begin{align}
{\rm Hess}\,f(S)[\xi] =S {\rm sym}( {\rm D} \nabla \bar{f}(S) [\xi] )S + {\rm sym} (\xi {\rm sym} (\nabla \bar{f}(S)) S ). \label{Hess}
\end{align}

The geodesic $\Gamma_{(S,\xi)}$ on the manifold ${\rm Sym}_+(r)$ going through a point $S\in {\rm Sym}_+(r)$ with a tangent vector $\xi\in T_S {\rm Sym}_+(r)$ is given by
\begin{align}
\Gamma_{(S,\xi)}(t) = \phi_{\exp (t \zeta)}(S) \label{geo1}
\end{align}
with $\xi = \zeta S + S \zeta$ for $\zeta\in T_{I_r} {\rm Sym}_+(r)$; i.e.,
the geodesic is the orbit of the one-parameter subgroup $\exp (t\zeta)$, where $\exp$ is the matrix exponential function.
The relation \eqref{geo1} follows from the fact that ${\rm Sym}_+(r)$ is a reductive homogeneous space.
For convenience, we prove it in Appendix \ref{ape1}. 
A detailed explanation of the expression \eqref{geo1} can be found in \cite{Gallier2016}.
To simplify \eqref{geo1}, we consider the geodesic going through the origin $I_r =\phi_{S^{-1/2}}(S) \in {\rm Sym}_+(r)$ because the Riemannian metric given by \eqref{metric} is invariant under the group action.
In this case, we get $\zeta=\frac{1}{2}\xi$ and $\Gamma_{(I_r,\xi)}(t) = \exp (t\xi)$.
Hence,
\begin{align*}
\Gamma_{(S,\xi)}(t) &= \phi_{S^{1/2}} (\Gamma_{(I_r, {\rm D}\phi_{S^{-1/2}}(\xi) )} (t))  \\
&= S^{\frac{1}{2}} \exp (t S^{-\frac{1}{2}} \xi S^{-\frac{1}{2}} ) S^{\frac{1}{2}}.
\end{align*}
Therefore, the exponential map on ${\rm Sym}_+(r)$ is given by
\begin{align}
{\rm Exp}_S (\xi) :=  \Gamma_{(S,\xi)}(1) 
= S^{\frac{1}{2}} \exp (S^{-\frac{1}{2}} \xi S^{-\frac{1}{2}} ) S^{\frac{1}{2}}. \label{8}
\end{align}
Since ${\rm Exp} :  T_S {\rm Sym}_+(r) \rightarrow {\rm Sym}_+(r)$ is a bijection \cite{Gallier2016},
 ${\rm Sym}_+(r)$ endowed with the Riemannian metric \eqref{metric} is geodesically complete in contrast to 
the case of \eqref{normal_metric}.

\subsection{Euclidean gradient of the objective function $J$}

Let $\bar{J}$ denote the extension of the objective function $J$ to the Euclidean space ${\bf R}^{r\times r}\times {\bf R}^{r\times m} \times {\bf R}^{p\times r}$.
Then, the Euclidean gradient of $\bar{J}$ is given by

\begin{align}
& \nabla \bar{J}(A_r,B_r,C_r) \nonumber \\
=& 2( -QP-Y^TX, QB_r+Y^TB, C_rP-CX). \label{16}
\end{align}

\noindent
Although a similar expression can be found in Theorem 3.3 in \cite{van2008h2} and Section 3.2 in \cite{wilson1970optimum},
we provide another proof in Appendix \ref{apeB} because some equations in the proof are needed for deriving the Riemannian Hessian of $J$ as shown in the next subsection.

\subsection{Geometry of Problem 1} \label{Sec:3.D}

We define the Riemannian metric of the manifold $M$ as
\begin{align}
& \langle (\xi_1,\eta_1,\zeta_1),(\xi_2,\eta_2,\zeta_2) \rangle_{(A_r,B_r,C_r)} \nonumber \\
:=& {\rm tr}( A_r^{-1} \xi_1 A_r^{-1} \xi_2  ) + {\rm tr} (\eta_1^T\eta_2) + {\rm tr}(\zeta_1^T \zeta_2) \label{Riemannian_metric}
\end{align}
for $(\xi_1,\eta_1,\zeta_1),(\xi_2,\eta_2,\zeta_2) \in T_{(A_r,B_r,C_r)} M$.
Then, it follows from \eqref{gradient} and \eqref{16} that
\begin{align}
 {\rm grad}\, J(A_r,B_r,C_r)  =& ( -2 A_r {\rm sym}(QP+Y^TX)A_r,  \label{grad_J} \\
& 2(QB_r+Y^TB), 2(C_rP-CX) ). \nonumber
\end{align}
Furthermore, from \eqref{Hess} and \eqref{16}, the Riemannian Hessian of $J$ at $(A_r,B_r,C_r)$ is given by
\begin{align}
&  {\rm Hess}\, J(A_r,B_r,C_r) [(A'_r,B'_r,C'_r)] \nonumber \\
=& ( -2A_r {\rm sym}( Q'P+QP'+Y'^TX+Y^TX') A_r \nonumber \\
&\, -2 {\rm sym} (A'_r {\rm sym} (QP+Y^TX) A_r), \label{Hess_J}\\
&\,\, 2(Q'B_r +QB'_r +Y'^TB), 2(C'_rP+C_rP'-CX') ), \nonumber
\end{align}
where $P'$ and $X'$ are the solutions to \eqref{10} and \eqref{11} in Appendix \ref{apeB}, respectively, and $Q'$ and $Y'$ are the solutions to
\begin{align}
& A_r Q' +Q' A_r+A'_r Q+Q A'_r- C'^T_r C_r-C^T_r C'_r =0, \label{31} \\
& AY'+Y'A_r+YA'_r+C^TC'_r =0. \label{32}
\end{align}
The equations \eqref{31} and \eqref{32} are obtained by differentiating \eqref{4} and \eqref{6}, respectively.
From \eqref{8}, we can define the exponential map on the manifold $M$ as
\begin{align}
& {\rm Exp}_{(A_r,B_r,C_r)}(\xi,\eta,\zeta) \nonumber \\
:=& (A_r^{\frac{1}{2}} \exp (A_r^{-\frac{1}{2}} \xi A_r^{-\frac{1}{2}} ) A_r^{\frac{1}{2}}, B_r+\eta,C_r+\zeta) \label{33}
\end{align}
for any $(\xi,\eta,\zeta)\in T_{(A_r,B_r,C_r)} M$;
 i.e., the manifold $M$ is geodesically complete.

\subsection{Trust-region method for Problem 1} \label{Sec:3.E}

This section gives the Riemannian trust-region method for solving Problem 1.
In \cite{absil2009optimization, absil2007trust}, the Riemannian trust-region method has been discussed in detail.

At each iterate $(A_r,B_r,C_r)$ in the Riemannian trust-region method on the manifold $M$, 
we evaluate the quadratic model $\hat{m}_{(A_r,B_r,C_r)}$ of the objective function $J$ within a trust-region:
\begin{align*}
&\quad \hat{m}_{(A_r,B_r,C_r)}(\xi,\eta,\zeta)  \\
=& J(A_r,B_r,C_r) + \langle {\rm grad}\,J(A_r,B_r,C_r), (\xi,\eta,\zeta) \rangle_{(A_r,B_r,C_r)} \\
&+\frac{1}{2}  \langle {\rm Hess}\, J(A_r,B_r,C_r)[(\xi,\eta,\zeta)], (\xi,\eta,\zeta) \rangle_{(A_r,B_r,C_r)}.
\end{align*}
A trust-region with a radius $\Delta>0$ at $(A_r,B_r,C_r)\in M$ is defined as a ball with center $0$ in $T_{(A_r,B_r,C_r)} M$.
Thus, the trust-region subproblem at $(A_r,B_r,C_r)\in M$ with a radius $\Delta$ is defined as a problem of minimizing $\hat{m}_{(A_r,B_r,C_r)}(\xi,\eta,\zeta)$ subject to $(\xi,\eta,\zeta)\in T_{(A_r,B_r,C_r)} M$, $||(\xi,\eta,\zeta)||_{(A_r,B_r,C_r)}:= \sqrt{ \langle (\xi,\eta,\zeta),(\xi,\eta,\zeta) \rangle_{(A_r,B_r,C_r)}} \leq \Delta$.
This subproblem can be solved by the truncated conjugate gradient method \cite{absil2009optimization}.
Then, we compute the ratio of the decreases in the objective function $J$ and the model $\hat{m}_{(A_r,B_r,C_r)}$ attained by the resulting $(\xi_*,\eta_*,\zeta_*)$ to decide
whether $(\xi_*,\eta_*,\zeta_*)$ should be accepted and whether the trust-region with the radius $\Delta$ is appropriate.
Algorithm \ref{algorithm} describes the process.
The constants $\frac{1}{4}$ and $\frac{3}{4}$ in the condition expressions in Algorithm \ref{algorithm} are commonly used in the trust-region method for a general unconstrained optimization problem.
These values ensure the convergence properties of the algorithm \cite{absil2009optimization, absil2007trust}.


\begin{algorithm}                      
\caption{Trust-region method for Problem 1.}    \label{algorithm}     
\label{alg1}                          
\begin{algorithmic}[1]
\STATE Choose an initial point $((A_{r})_0,(B_{r})_0,(C_{r})_0) \in M$ and parameters $\bar{\Delta}>0$, $\Delta_0\in (0,\bar{\Delta})$, $\rho'\in [0,\frac{1}{4})$.
\FOR{$k=0,1,2,\ldots$ }
\STATE Solve the following trust-region subproblem for $(\xi,\eta,\zeta)$ to obtain $(\xi_k,\eta_k,\zeta_k)\in T_{(A_r,B_r,C_r)} M$:
\begin{align*}
&{\rm minimize}\quad \hat{m}_{((A_r)_k,(B_r)_k,(C_r)_k)}(\xi,\eta,\zeta) \\
&{\rm subject\, to}\quad ||(\xi,\eta,\zeta)||_{((A_r)_k,(B_r)_k,(C_r)_k)} \leq \Delta_k, \\
&{\rm where}\quad \hat{m}_k(\xi,\eta,\zeta):=\hat{m}_{((A_r)_k,(B_r)_k,(C_r)_k)}(\xi,\eta,\zeta), \\
&\quad\quad\quad  (\xi,\eta,\zeta)\in T_{((A_r)_k,(B_r)_k,(C_r)_k)}M.
\end{align*}
\STATE Evaluate
\begin{align*}
\rho_k := \frac{ J({\rm Exp}_{k}(0,0,0)) -J({\rm Exp}_{k}(\xi_k,\eta_k,\zeta_k))}{ \hat{m}_{k}(0,0,0)- \hat{m}_{k} (\xi_k,\eta_k,\zeta_k)}
\end{align*}
\STATE with
${\rm Exp}_k(\xi,\eta,\zeta):= {\rm Exp}_{((A_r)_k,(B_r)_k,(C_r)_k)}(\xi,\eta,\zeta)$.
\IF {$\rho_k<\frac{1}{4}$}
\STATE 
$\Delta_{k+1}=\frac{1}{4}\Delta_k$.
\ELSIF {$\rho_k>\frac{3}{4}$ and $||(\xi_k,\eta_k,\zeta_k)||_{((A_r)_k,(B_r)_k,(C_r)_k)} = \Delta_k$}
\STATE
$\Delta_{k+1} = \min (2\Delta_k,\bar{\Delta})$.
\ELSE 
\STATE
$\Delta_{k+1} = \Delta_k$.
\ENDIF
\IF {$\rho_k>\rho'$}
\STATE
$((A_r)_{k+1},(B_r)_{k+1},(C_r)_{k+1}) = {\rm Exp}_k(\xi_k,\eta_k,\zeta_k)$.
\ELSE
\STATE
{\small
$((A_r)_{k+1},(B_r)_{k+1},(C_r)_{k+1}) = ((A_r)_{k},(B_r)_{k},(C_r)_{k})$.}
\ENDIF
\ENDFOR
\end{algorithmic}
\end{algorithm}

\begin{remark} \label{remark3}
The most computational task to perform Algorithm 1 is to solve Eqs. \eqref{3}--\eqref{6} iteratively.
Although some algorithms to solve these equations have been studied in some literatures \cite{antoulas2005approximation, benner2003state, damm2008direct},
we need to develop a more effective method for solving large-scale model reduction problems by Algorithm 1.
\end{remark}

\section{Optimization algorithm for solving Problem 2} \label{sec_Pro2}
This section develops an optimization algorithm for solving Problem 2.

As with Problem 1, to derive the Riemannian gradient and Hessian of the objective function $\tilde{J}$,
we calculate the Euclidean gradient $\nabla \bar{\tilde{J}}$, where $\bar{\tilde{J}}$ is the extension of $\tilde{J}$ to the ambient space ${\bf R}^{r\times r}\times {\bf R}^{r\times m}$.
Since $C=B^T$ and $C_r=B_r^T$, it follows from \eqref{3}--\eqref{6} that
$P=Q$ and $X=-Y$.
Thus, in Appendix \ref{apeB}, by replacing $C$, $C_r$, and $C'_r$ with $B$, $B_r$, and $B'_r$, respectively,
we obtain
\begin{align*}
\nabla \bar{\tilde{J}}(A_r,B_r) 
= ( -2(P^2-X^TX), 4PB_r-4X^TB).
\end{align*}
Hence, if we consider the counterpart of the Riemannian metric \eqref{Riemannian_metric} for the manifold $\tilde{M}$ as
\begin{align*}
 \langle (\xi_1,\eta_1),(\xi_2,\eta_2) \rangle_{(A_r,B_r)} = {\rm tr}( A_r^{-1} \xi_1 A_r^{-1} \xi_2  ) + {\rm tr} (\eta_1^T\eta_2),
\end{align*}
the Riemannian gradient and Hessian of $\tilde{J}$ are given by 
\begin{align*}
& {\rm grad}\, \tilde{J}(A_r,B_r) = (-2A_r{\rm sym}(P^2 -X^T X) A_r, \\
&\quad \quad\quad\quad\quad\quad\quad \quad 4PB_r-4X^TB ), \\
&  {\rm Hess}\, \tilde{J}(A_r,B_r) [ (A'_r,B'_r)]  \\
=& ( -2A_r {\rm sym}( P'P+PP'-X'^TX-X^TX') A_r  \\
&\, -2 {\rm sym} (A'_r {\rm sym} (P^2-X^TX) A_r), \\
&\,\, 4(P'B_r+PB'_r)  -4 X'^TB ), 
\end{align*}
respectively.
Here,
$P$, $X$, $P'$, and $X'$ are the solutions to \eqref{3}, \eqref{5}, \eqref{10}, and \eqref{11}, respectively.
The exponential map on the manifold $\tilde{M}$ is, of course, given by
\begin{align*}
 {\rm Exp}_{(A_r,B_r)}(\xi,\eta) =(A_r^{\frac{1}{2}} \exp (A_r^{-\frac{1}{2}} \xi A_r^{-\frac{1}{2}} ) A_r^{\frac{1}{2}}, B_r+\eta).
\end{align*}
Similarly to Problem 1, we can solve Problem 2 by using a modified algorithm of Algorithm \ref{algorithm}. The reduced system constructed by the solution is also a stable gradient system.
Note that in contrast to Problem 1, we do not calculate $Q$, $Y$, $Q'$, and $Y'$; i.e., we only need to calculate $P$, $X$, $P'$, and $X'$ for solving Problem 2 by the trust-region method.
This improves computational efficiency.

\begin{remark}
If $m=p$ and $B^T=C$, we can also regard the system \eqref{1} as a port-Hamiltonian system \cite{ionescu2013moment, van2012l2}.
Since a port-Hamiltonian system is passive, the reduced system constructed by the solution to Problem 2 is also passive \cite{van2012l2}.
\end{remark}

\section{Comparison between Problem 1 and Problem 3} \label{secnew_5}
In this section, we compare the reduced systems obtained by solving Problems 1 and 3 and give a simple example which shows that they do not necessarily coincide with each other.

For $J$ in \eqref{eq_J}, which is the $H^2$ norm of the error system, let $J_1 := J$ and $J_3(U) := J(U^TAU, U^TB, CU)$.
Then, we have
\begin{equation*}
{\rm grad}\, J_1(A_r, B_r, C_r)=(A_r{\rm sym}(\nabla_{A_r} \bar{J})A_r, \nabla_{B_r} \bar{J}, \nabla_{C_r} \bar{J})
\end{equation*}
and
\begin{equation*}
{\rm grad}\, J_3(U) = \nabla \bar{J_3}(U) - U {\rm sym}(U^T \nabla \bar{J_3}(U)),
\end{equation*}
where
\begin{align*}
\nabla \bar{J_3}(U) =& 2 A U {\rm sym}(\nabla_{A_r} \bar{J}(U^TAU, U^TB, CU)) \\
&+ B(\nabla_{B_r} \bar{J}(U^TAU, U^TB, CU))^T \\
&+ C^T\nabla_{C_r} \bar{J}(U^TAU, U^TB, CU).
\end{align*}
Note that we have used $A^T=A$ and that $\nabla_{A_r} \bar{J}$ denotes the $A_r$-component of $\nabla \bar{J}$.
The expression of ${\rm grad}\, J_1(A_r,B_r,C_r)$ is from \eqref{16} and \eqref{grad_J}, and ${\rm grad}\, J_3(U)$ can be found in \cite{sato2015riemannian, yan1999approximate}.

Even if ${\rm grad}\, J_1(A_r, B_r, C_r) = 0$ for some $(A_r, B_r, C_r)$, there does not in general exist $U$ such that
\begin{equation}
\label{EqABC}
A_r = U^T A U,\ B_r = U^TB,\ C_r = CU,
\end{equation}
and 
${\rm grad}\,J_3(U) =0$.
Conversely, ${\rm grad}\, J_3(U) =0$ does not yield ${\rm grad}\, J_1(U^TAU, U^TB, CU)=0$ either.

In order to see this clearly from a simple example, we consider in the remainder of this section the system \eqref{1} with $n=2$ and $m=p=1$ and assume that the dimension of the reduced model is $r=1$.
Furthermore, we suppose $A=\begin{pmatrix}2 & 0 \\ 0 & 1\end{pmatrix}$, $B = \begin{pmatrix}-1 \\ 1\end{pmatrix}$, and $C = \begin{pmatrix} 1 & 1 \end{pmatrix}$.

For Problem 1, we can obtain $P = B_r^2/2A_r$, $Q = C_r^2/2A_r$, $X = \begin{pmatrix}-B_r/(A_r+2) & B_r/(A_r+1)\end{pmatrix}^T$, and $Y = -\begin{pmatrix}C_r/(A_r+2) & C_r/(A_r+1)\end{pmatrix}^T$ by \eqref{3}--\eqref{6}.
Then, a simple analysis implies that ${\rm grad}\, J_1(A_r, B_r, C_r) = 0$ is equivalent to
\begin{equation*}
B_r C_r = 0 \quad \text{or} \quad A_r=-\frac{1}{2}+\frac{\sqrt{33}}{6},\ B_rC_r=6-\sqrt{33}.
\end{equation*}
The objective function at these infinite critical points are evaluated as
\begin{align*}
J(A_r, B_r, C_r) = 1/12 = 0.0833
\end{align*}
for any $(A_r, B_r, C_r)$ with $A_r>0$ and $B_rC_r=0$, and
\begin{equation*}
J(A_r, B_r, C_r) = (569-99\sqrt{33})/24 = 0.0120
\end{equation*}
for $A_r=-1/2+\sqrt{33}/6$ and for any $(B_r, C_r)$ with $B_rC_r=6-\sqrt{33}$, which implies that the minimum value of $J$ attained by solving Problem 1 is $0.0120$.

For Problem 3, let $U = \begin{pmatrix} u_1 & u_2 \end{pmatrix}^T \in {\rm St}\,(1,2)$.
This means that $U$ is in the unit $2$-sphere, that is, $u_1^2+u_2^2=1$.
Then, we have $P = (u_1-u_2)^2/2(1+u_1^2)$, $Q = (u_1+u_2)^2/2(1+u_1^2)$,\\
$X =\begin{pmatrix} (u_1-u_2)/(u_1^2+3) & - (u_1-u_2)/(u_1^2+2)\end{pmatrix}^T$,
and
$Y = -\begin{pmatrix}(u_1+u_2)/(u_1^2+3) & (u_1+u_2)/(u_1^2+2)\end{pmatrix}^T$ in a similar manner to that in Problem 1.
A straightforward but tedious calculation shows that ${\rm grad}\, J_3(U) = 0$ holds if and only if
\begin{equation}
\label{Ueq2}
(u_1, u_2) = (\pm 1, 0),\ (0, \pm 1)
\end{equation}
or
\begin{equation}
\label{Ueq1}
u_1 = \pm 0.5642 \quad \text{and} \quad u_2 = \pm \sqrt{1-u_1^2}= \pm 0.8256,
\end{equation}
where $u_1=\pm 0.5642$ are the real solutions to the equation
$4u_1^{12}+48u_1^{10}+215u_1^8+478u_1^6+515u_1^4+132u_1^2-112=0$.
Therefore, there are only $8$ finite discrete critical points of $J_3$ in contrast to Problem 1.
The resultant reduced system matrices are then computed by \eqref{EqABC}.
Eq.~\eqref{Ueq2} yields $(A_r, B_r, C_r) = (2, \mp 1, \pm 1), (1, \pm 1, \pm 1)$, where $B_r C_r = \pm 1$.
In contrast, for \eqref{Ueq1} we  have $A_r = 1.318$ and $(B_r, C_r) = (\pm 0.2614, \pm 1.390)$, $(\pm 1.390, \pm 0.2614)$, where $B_rC_r=0.3633$.
Meanwhile, the result for Problem 1 yields $B_rC_r=0$ or $B_rC_r=6-\sqrt{33}=0.2554$.
Therefore, we can conclude that the reduced systems obtained by the two problems do not coincide with each other in general.
Furthermore, we have $J(A_r,B_r,C_r)=0.0389$ for all $(A_r, B_r, C_r)$ obtained by \eqref{Ueq1}, $J(A_r, B_r, C_r) =  1/2 = 0.5$ for $(u_1, u_2) = (\pm 1, 0)$, and $J(A_r, B_r, C_r) = 1/4 = 0.25$ for $(u_1, u_2) = (0, \pm 1)$,
all of which are worse than the results in Problem 1.

From these observations, we can conclude that the solutions to Problems 1 and 3 are not necessarily unique nor the solution sets of both problems do not contain each other.
Also, the attained optimal values do not coincide with each other.

\section{Numerical experiments} \label{sec4}
This section illustrates that the proposed reduction method preserves the structure of the system \eqref{1} although the balanced truncation method does not preserve it.
Furthermore, it is shown that the value of the objective function in the case of the proposed reduction method becomes smaller  than that in the case of the reduction method proposed in \cite{sato2015riemannian} even if we choose an initial point in Algorithm 1 as a local optimal solution to Problem 3.
This means that the stationary points of Problems 1 and 3 do not coincide.
To perform them, we have used Manopt \cite{boumal2014manopt}, which is a MATLAB toolbox for optimization on manifold.

We consider a reduction of the system \eqref{1} with $n=5$ and $m=p=2$ to the system \eqref{2} with $r=3$.
Here, the system matrices $A$, $B$, and $C$ are given by
\begin{align*}
A &:= \begin{pmatrix}
3 & -1 & 1 & 1 & -1 \\
-1 & 2 & 0 & 0 & 2 \\
1 & 0 & 2 & 1 & 1 \\
1 & 0 & 1 & 3 & 0 \\
-1 & 2 & 1 & 0 & 4
\end{pmatrix},
B := \begin{pmatrix}
0 & 1 \\
1 & 0 \\
-1 & 1 \\
1 & 0 \\
0 & 1
\end{pmatrix},\\
C &:= \begin{pmatrix}
1 & 0 & 0 & 0 & 0 \\
0 & 0 & 1 & 0 & 1
\end{pmatrix}.
\end{align*}
That is, $(A,B,C)\in {\rm Sym}_+(5)\times {\bf R}^{5\times 2} \times {\bf R}^{2\times 5}$.

The balanced truncation method, which is the most popular model reduction method \cite{antoulas2005approximation}, gave the reduced matrix $A_r^{\rm BT}$ as
\begin{align*}
A_r^{\rm BT} &= \begin{pmatrix}
2.8944 & -0.0422 & -1.4729 \\
-0.0318 & 1.0470 & -0.2615 \\
-1.1764 & -0.2355 & 4.1898
\end{pmatrix}.
\end{align*}
Thus, $A_r^{\rm BT} \not\in {\rm Sym}_+(3)$; i.e., the balanced truncation method did not preserve the original model structure.
Furthermore, we obtained $||G-G_r||_{H^2}=0.0157$.

The reduction method which was briefly explained in Remark 1 in \cite{sato2015riemannian}, gave the orthogonal matrix
\begin{align}
U = \begin{pmatrix}
0.8906 & 0.1189 & -0.1025 \\
-0.1117 & 0.7216 & 0.0373 \\
-0.0650 & -0.1558 & 0.8994 \\
-0.2144 & 0.6138 & 0.0302 \\
0.3798 & 0.2532 & 0.4223
\end{pmatrix}, \label{U}
\end{align}
and then
\begin{align*}
U^TAU &= 
\begin{pmatrix}
1.9613 & 0.0507 & 0.7510 \\
0.0507 & 2.8566 & 1.6666 \\
0.7510 & 1.6666 & 3.1486
\end{pmatrix}.
\end{align*}
Thus, $U^TAU \in {\rm Sym}_+(3)$; i.e., this method preserved the original model structure.
Furthermore, we obtained $||G-G_r||_{H^2}=0.0217$.
Note that, in this result, the norm of the gradient of the objective function was approximately equal to $7.493\times 10^{-7}$; i.e.,
we can expect that a local optimal solution to Problem 3 was obtained.

The proposed algorithm gave the reduced matrix $A_r$, $B_r$, and $C_r$ as follows:
\begin{align*}
A_r &= \begin{pmatrix}
1.8965  &  0.0237 &  0.7778 \\
0.0237 &   3.1554  &  1.8009 \\
0.7778 &   1.8009 &   3.1784
\end{pmatrix}, \\
B_r & = \begin{pmatrix}
  -0.2677 &  1.1820 \\
    1.5124  &  0.2049 \\
  -0.7759  &  1.2155
\end{pmatrix}, \\
C_r &= \begin{pmatrix}
  0.8726 &  0.1503 &  -0.0630 \\
  0.3321 &  0.0680 &  1.3121
\end{pmatrix}.
\end{align*}
Thus, $A_r\in {\rm Sym}_+(3)$; i.e., the reduced system had the same structure with the original system.
Here, we chose an initial point $((A_r)_0, (B_r)_0, (C_r)_0)$  in Algorithm \ref{algorithm} as $(U^TAU, U^TB, CU)$,
where $U$ is defined by \eqref{U}.
Furthermore, we obtained $||G-G_r||_{H^2}=0.0156$.
Hence, the value of the objective function attained by the proposed algorithm was smaller than those by the balanced truncation method and the method in \cite{sato2015riemannian}.
This means that the stationary points of Problems 1 and 3 do not coincide.

To verify the effectiveness of the proposed algorithm for medium-scale systems,
we also randomly created matrices $A$, $B$, and $C$ of larger size.
Table \ref{table} shows the values of the relative $H^2$ error in the case of $A\in {\rm Sym}_+(300)$, $B\in {\bf R}^{300\times 3}$, and $C\in {\bf R}^{2\times 300}$, respectively.
For all $r$, the relative $H^2$ errors in the proposed method were smaller than those of the balanced truncation method.
Furthermore, the reduced models by the balanced truncation method did not have the original symmetric structure while the proposed method had.
Moreover, for all $r$, 
 the proposed method was better than the method in \cite{sato2015riemannian}.
Here, we note that for each $r$, an initial point $((A_r)_0, (B_r)_0, (C_r)_0)$  in Algorithm \ref{algorithm} to solve Problem 1 was chosen as $(U^TAU, U^TB, CU)$,
where $U$ is a local optimal solution to Problem 3.
Thus, Table \ref{table} also shows that the stationary points of Problems 1 and 3 do not coincide.

\begin{table}[h]
\caption{The comparison of the relative $H^2$ error $\frac{||G-G_r||_{H^2}}{||G||_{H^2}}$.} \label{table}
  \begin{center}
    \begin{tabular}{|c|c|c|c|c|} \hline
         $r$&   $6 $   &  $8$ & $10$ &  $12$  \\ \hline 
      Balanced truncation                              & 0.0141  & 0.0120 & 0.0103 & 0.0088 \\\hline 
	The method in \cite{sato2015riemannian} & 0.0297 & 0.0299 & 0.0294 & 0.0317 \\\hline 
The proposed method                                 & 0.0112 & 0.0089& 0.0042 & 0.0020 \\\hline 
    \end{tabular}
  \end{center}
\end{table}

\begin{remark}
As mentioned in Remark \ref{remark3}, in order to solve large-scale model reduction problems by Algorithm 1,
a long computational time is needed.
On the other hand, a computational time for performing the balanced truncation method is less than it.
Furthermore, the balanced truncation method gives upper bounds of the $H^2$ and $H^{\infty}$ error norms \cite{antoulas2005approximation}.
From these facts, we suggest that we use the balanced truncation method for determining the possible largest reduced dimension $r$ for performing Algorithm 1 by observing the $H^2$ and $H^{\infty}$ error norms.
Then, we can choose an actual $r$ to perform Algorithm 1 as a smaller value than the possible largest dimension.
\end{remark}

\section{Conclusion} \label{sec5}

We have studied the stability and symmetry preserving $H^2$ optimal model reduction problem on the product manifold of the manifold of the 
symmetric positive definite matrices and two Euclidean spaces.
To solve the problem by using the trust-region method, we have derived the Riemannian gradient and Riemannian Hessian.
Furthermore, it has been shown that if we restrict our systems to gradient systems, the gradient and Hessian can be obtained more efficiently.
By a simple example, we have proved that 
the solutions to our problem and the problem in \cite{sato2015riemannian} are not unique and the solution sets of both problems do not contain each other in general.
Also, it has been revealed that the attained optimal values do not coincide.
Numerical experiments have illustrated that although the balanced truncation does not preserve the original symmetric structure of the system, the proposed method preserves the structure.
Furthermore,  it has been demonstrated that
the proposed method is better than our method in \cite{sato2015riemannian}, and also usually better than the balanced truncation method, in the sense of the  $H^2$ error norm between the transfer functions of the original and reduced systems.


%

\appendix

\subsection{Proof of the fact that ${\rm Sym}_+(r)$ is a reductive homogeneous space} \label{ape1}

To prove that ${\rm Sym}_+(r)$ is a reductive homogeneous space, we first note that
there is a natural bijection
\begin{align}
{\rm Sym}_+(r) \cong GL(r)/O(r). \label{doukei}
\end{align}
To see this, let $\phi_g$ be
$GL(r)$ action on the manifold ${\rm Sym}_{+}(r)$; i.e.,
$\phi_g(S) =gSg^T,\quad g\in GL(r), S\in {\rm Sym}_{+}(r)$.
The action $\phi_g$ is transitive; i.e., for any $S_1$, $S_2\in {\rm Sym}_+(r)$, there exists $g\in GL(r)$ such that
$\phi_g(S_1)=S_2$.
Thus, the manifold ${\rm Sym}_+(r)$ consists of a single orbit; i.e., ${\rm Sym}_+(r)$ is a homogeneous space of $GL(r)$.
The action $\phi_g$ has the isotropy subgroup of the orthogonal group $O(r)$ at  $I_r\in {\rm Sym}_+(r)$ because
$O(r) = \{ g\in GL(r)\,|\, \phi_g (I_r) = I_r \}$.
In general, if an action of a group on a set is transitive, the set is isomorphic to a quotient of the group by its isotropy subgroup \cite{Gallier2016}.
Hence, \eqref{doukei} holds.
From the identification \eqref{doukei}, we can show that the quotient $GL(r)/O(r)$ is reductive; i.e.,
$T_{I_r} GL(r) \cong T_{I_r} {\rm Sym}_+(r) \oplus T_{I_r} O(r)$ 
and 
$O\xi O^{-1} \in T_{I_r} {\rm Sym}_+(r)$ 
for $\xi\in T_{I_r} {\rm Sym}_+(r)$ and $O\in O(r)$.
In fact, 
these follow from
\begin{align*}
T_{I_r} GL(r) &\cong {\bf R}^{r\times r} \cong {\rm Sym}(r) \oplus {\rm Skew}(r), \\
{\rm Sym}(r) & \cong T_{I_r} {\rm Sym}_+(r), \,\,
{\rm Skew} (r) \cong T_{I_r} O(r),
\end{align*}
and $O^{-1}=O^T$.

\subsection{Proof of \eqref{16}} \label{apeB}

The directional derivative of $\bar{J}$ at $(A_r,B_r,C_r)$ in the direction $(A'_r,B'_r,C'_r)$ can be calculated as
\begin{align}
&{\rm D}\bar{J}(A_r,B_r,C_r)[(A'_r,B'_r,C'_r)]  \nonumber \\
=& 2{\rm tr} (C'_r (PC_r^T -X^TC^T)) - 2{\rm tr} (C^TC_rX'^T) +{\rm tr} (C_rP'C_r^T), \label{9}
\end{align}
where $P'$ and $X'$ are also the directional derivative of $P$ and $X$ at $(A_r,B_r,C_r)$ in the direction $(A'_r,B'_r,C'_r)$, respectively.
Differentiating \eqref{3} and \eqref{5}, we obtain 
\begin{align}
& A_r P'+P' A_r +A_r' P+PA'_r-B_r'B_r^T-B_rB'^{T}_r =0, \label{10} \\
& AX'+X'A_r+XA_r'-BB_r'^T =0. \label{11}
\end{align}
Eqs.\,\eqref{4} and \eqref{10} yield that
\begin{align}
{\rm tr}(C_r^TC_rP') = -2 {\rm tr}(A_r'^TQP)+2{\rm tr}(B_r'^TQB_r), \label{12}
\end{align}
and \eqref{6} and \eqref{11} imply that
\begin{align}
{\rm tr}(-C^TC_rX'^T) = {\rm tr}( (-X A_r'^T+BB_r'^T)^TY). \label{13}
\end{align}
By substituting \eqref{12} and \eqref{13} into \eqref{9}, we have
\begin{align}
&{\rm D}\bar{J}(A_r,B_r,C_r)[(A'_r,B'_r,C'_r)]  \nonumber \\
=& 2{\rm tr} ( A_r'^T (-QP-Y^TX) ) + 2{\rm tr} (B_r'^T (QB_r+Y^TB) ) \nonumber \\
&+2{\rm tr} (C_r'^T(C_rP-CX) ). \label{14}
\end{align}
Since the Euclidean gradient $\nabla \bar{J}(A_r,B_r,C_r)$ satisfies
\begin{align*}
& {\rm D} \bar{J}(A_r,B_r,C_r)[(A'_r,B'_r,C'_r)]  \nonumber \\
= & {\rm tr} (A_r'^T \nabla_{A_r} \bar{J}(A_r,B_r,C_r)) + {\rm tr} (B_r'^T \nabla_{B_r} \bar{J}(A_r,B_r,C_r))\\
&+ {\rm tr} (C_r'^T \nabla_{C_r} \bar{J}(A_r,B_r,C_r)),
\end{align*}
\eqref{14} implies \eqref{16}.

\section*{Acknowledgment}

This study was supported in part by JSPS KAKENHI Grant Number JP16K17647.
The authors would like to thank the anonymous reviewers
for their valuable comments that helped improve the paper
significantly.

\ifCLASSOPTIONcaptionsoff
  \newpage
\fi

\end{document}